\documentclass[12pt]{article}

\usepackage[english]{babel}
\usepackage{amssymb,amsmath}
\usepackage{hyperref}

\newcommand{\footremember}[2]{%
    \footnote{#2}
    \newcounter{#1}
    \setcounter{#1}{\value{footnote}}%
}
\newcommand{\footrecall}[1]{%
    \footnotemark[\value{#1}]%
}

\newtheorem{theorem}{Theorem }[section]
\newtheorem{lemma}[theorem]{Lemma}

\newtheorem{proposition}[theorem]{Proposition}

\newcommand{\proof}{\noindent\textbf{Proof. }}
\newcommand{\qed}{\hspace*{\fill}$\Box$}

\providecommand{\keywords}[1]{\textbf{Keywords---} #1}

\newcommand{\cB}{\mathcal{B}}

\title{Non-existence of  partial difference sets of order $8p^3$ in Abelian groups}

\author{Stefaan De Winter \footremember{MTU}{Michigan Technological University} \footnote{sgdewint@mtu.edu} \and Zeying Wang \footrecall{MTU} \footnote{zeying@mtu.edu}}

\date{}

\begin{document}

\maketitle 

\abstract{In this paper we prove non-existence of  nontrivial  partial difference sets in Abelian groups of order $8p^3$, where  $p\ge 3$ is  a prime number.}


\section{Introduction}

Let $G$ be a finite Abelian group of order $v$, and let $D\subseteq G$ be a subset of size $k$. We say $D$ is a $(v,k,\lambda,\mu)$-{\it partial difference set} (PDS) in $G$ if the expressions $gh^{-1}$, $g\in D$, $h\in D$, $g\neq h$, represent each non-identity element in $D$ exactly $\lambda$ times, and each non-identity element of $G$ not in $D$ exactly $\mu$ times. Further assume that $D^{(-1)}=D$ (where $D^{(s)}=\{g^s:g\in D\}$ ) and $e\notin D$, where $e$ is the identity of $G$, then $D$ is called a {\it regular} partial difference set.  A regular PDS is called {\it trivial} if $D\cup\{e\}$ or $G\setminus D$ is a subgroup of $G$. The condition that $D$ be regular is not a very restrictive one, as  $D^{(-1)}=D$ is automatically fulfilled whenever $\lambda\neq\mu$, and $D$ is a PDS if and only if $D\cup\{e\}$ is a PDS. The importance of regular PDS lies in the fact that they are equivalent to strongly regular Cayley graphs. A detailed study of PDS was started by Ma in \cite{MA84}. By now there is a rich literature on PDS, both with a focus on existence conditions and classification (see for example \cite{Arasu}, \cite{SDWZW}, \cite{Leung2}, \cite{MA97},\dots), as well as with a focus on constructions (see for example \cite{Davis}, \cite{Hou}, \cite{Pol}, \cite{Pol2}, \dots). 
\medskip

When the group $G$ is Abelian only a limited number of examples of regular PDS are known: (negative) Latin square type PDS, reversible difference sets, Paley type PDS, PCP type PDS, and projective two-weight sets. However, all these examples have many connections to other areas in combinatorics, and so it is natural to look for other examples. On the other hand, not that many classification results are known, leaving an abundance of parameter sets for which no PDS is known, but for which existence has not been excluded. In recent work the authors proved non-existence of PDS in 18 specific cases \cite{SDWEKZW,SDWZW_216}, finalizing the classification of parameters for which there exists a regular PDS with $k\leq 100$ in an Abelian group. This classification was started by Ma \cite{MA94b}, and had been open for almost 20 years. Building on some of the techniques developed in \cite{SDWEKZW} the authors obtained in \cite{SDWZW} a complete classification of regular PDS in Abelian groups of order $4p^2$, $p\geq3$ a prime. Here we continue the classification of PDS in Abelian groups, based on the order of the group.  Our main result is that in Abelian groups of order $8p^3$, $p\geq3$ a prime, only trivial PDS exist. Surprisingly the technique used is very different from that used in the classification for groups of order $4p^2$. The main reason for this is that in the $4p^2$ case the number of hypothetical parameter sets for each $p$ grows rapidly as $p$ grows, whereas in the $8p^3$ case this is not the case. When searching for parameters that survive the basic integrality and divisibility conditions we noticed that for most primes $p$ only two parameter sets appeared, and for the remaining primes exactly 4 sets appeared. This explains why we use an approach focussed on restricting parameters in this paper (Section 2), rather than the constructive approach taken in \cite{SDWZW}. Equally surprising to the authors was the fact that some of the seemingly ad hoc methods used to prove nonexistence of PDS in Abelian groups of order $2^3\cdot 3^3$ remain largely valid in the general case (Section 3). The results of this paper seem to further confirm that PDS in Abelian groups are actually rare objects.

\medskip

We end this section by listing a few useful results which we will need in the rest of the paper. The first three results can be found in Ma (\cite{MA94b}, \cite{MA97}), and the fourth is a recent local multiplier theorem from \cite{SDWEKZW}. Given a $(v,k,\lambda,\mu)$-PDS, one defines the following two parameters: $\beta:=\lambda-\mu$, and $\Delta:=(\lambda-\mu)^2+4(k-\mu)$.

\begin{proposition}\label{Ma2}
No non-trivial PDS exists in
\begin{itemize}
\item an Abelian group $G$ with a cyclic Sylow-$p$-subgroup and $o(G)\neq p$;
\item an Abelian group $G$ with a Sylow-$p$-subgroup isomorphic to $\mathbb{Z}_{p^s}\times\mathbb{Z}_{p^t}$ where $s\neq t$.
\end{itemize}
\end{proposition}

\begin{proposition}\label{Ma1}
Let D be a nontrivial regular $(v,\,k,\,\lambda,\,\mu)$-PDS in an Abelian group G. Suppose $\Delta=(\lambda-\mu)^2+4(k-\mu)$ is a perfect square.  If $N$ is a subgroup of $G$ such that $\gcd(\left|N\right|, \left|G\right|/\left|N\right|)=1$ and $\left|G\right|/\left|N\right|$ is odd, then $D_1=D\cap N$ is a (not necessarily non-trivial)  regular $(v_1,\,k_1,\,\lambda_1,\,\mu_1)$-PDS with
 $$\left|D_1\right|=\frac{1}{2}\left[ \left|N\right| +\beta_1\pm \sqrt{(\left|N\right|+\beta_1)^2-(\Delta_1-\beta_1^2)(\left|N\right|-1)}  \right].$$
Here $\Delta_1=\pi^2$ with $\pi=\gcd(\left|N\right|,\sqrt{\Delta})$ and $\beta_1=\beta-2\theta\pi$ where $\beta=\lambda-\mu$ and $\theta$ is the integer satisfying $(2\theta-1)\pi\leq\beta<(2\theta+1)\pi$.
\end{proposition}

\begin{proposition}\label{Ma97}
 Suppose there exists a regular ($v$, $k$, $\lambda$, $\mu$)-PDS $D$ in a group $G$. 
\begin{itemize}
\item[\rm{ (a)}] If $D \neq \emptyset$ and $D\neq G\setminus\{e\}$, then $0 \le \lambda \le k-1$ and $0\le \mu \le  k-1$.
\item[{\rm (b)}] The parameters $\beta$ and $\Delta$ have the same parity.
\item[{\rm (c)}] The PDS $D$  is nontrivial if and only if $-\sqrt{\Delta }<\beta <\sqrt{\Delta}-2$.  Also, if $D\neq G\setminus\{e\}$, then $D$ is nontrivial if and only if $1 \le \mu \le k-1$.
\item[{\rm (d)}] If $\Delta$  is not a square, then ($v$, $k$, $\lambda$, $\mu$)=($4t+1$, $2t$, $t-1$, $t$) for some positive integer $t$; furthermore, if $G$ is Abelian, then $v=p^{2s+1}$ for some prime $p\equiv 1 \pmod 4$.
\item[{\rm (e)}] If G is Abelian and $D\neq \emptyset$ and  $D\neq G\setminus\{e\}$, then $v^2\equiv (2k-\beta)^2\equiv 0 \pmod \Delta$; furthermore, if $D$ is nontrivial, then $v$, $\Delta$, and $v^2/ \Delta$ have the same prime divisors.
\item[{\rm (f)}] The set $(G\setminus  D)\setminus \{e\}$ is a PDS with parameters $(v', \, k', \, \lambda',\, \mu')$=$(v, \, v-k-1, \, v-2k-2+\mu, \, v-2k+\lambda)$ called the complement of $D$.
\item[{\rm (g)}] If $D$ is nontrivial, then there exists a nontrivial regular $(v,k^+,\lambda^+,\mu^+)$-PDS $D^+$ (in an Abelian group of order $v$) with $\Delta^+=(\lambda^+-\mu^+)^2+4(k^+-\mu^+)=v^2/\Delta$. (The PDS $D^+$ is called the {\it dual} of $D$.)
\end{itemize}
\end{proposition}

\begin{proposition}\label{lmt}{\rm [LMT]}
Let D be a regular $(v,k,\lambda,\mu)$-PDS in an Abelian group $G$. Furthermore assume $\Delta=(\lambda-\mu)^2+4(k-\mu)$ is a perfect square.  Then $g\in G$ belongs to $D$ if and only if $g^s$ belongs to $D$ for all $s$ coprime with $o(g)$, the order of $g$.
\end{proposition}

Throughout the rest of the paper we will assume that $p$ is a prime with $p\geq5$. The reason for this is twofold: the case $p=3$ already was dealt with in \cite{SDWZW_216}, and some of our arguments are only valid if $p\geq5$. Furthermore, we will always assume that $D$ is a nontrivial regular PDS.

\section{Parameter Restrictions}

\subsection{Restriction on $k$}

By part (f) of Proposition \ref{Ma97} we may assume that $k\leq v/2$. By part (g) of Proposition \ref{Ma97}, we may assume that $\Delta\leq v$. Before proving a better upper bound on $k$ we will show that under our assumptions $\Delta$ can only take two values. 

\begin{lemma}\label{Delta}
If $D$ is a regular ($v$, $k$, $\lambda$, $\mu$)-PDS in an Abelian group with $v=8p^3$, $k\leq v/2$ and $\Delta\leq v$ then either $\Delta=4p^2$ or $\Delta=16p^2$.
\end{lemma}

\proof  By Ma's Proposition \ref{Ma97} (d) and (e),  $\Delta$ must be a square, and $v$, $\Delta$, $v^2/ \Delta$ have the same prime divisors. It follows that $\Delta=4p^2$ or $\Delta=16p^2$. \qed

\medskip

Since the Cayley graph of a $(v,\,k,\,\lambda, \,\mu)$--PDS is a $(v,\,k,\,\lambda, \,\mu)$-strongly regular graph, we have 
\begin{equation}\label{srg_cond_1}
k(k-\lambda-1)=\mu(v-k-1).
\end{equation}

By substituting $\frac{k(k-\lambda-1)}{v-k-1}$ for $\mu$ into $\Delta=(\lambda-\mu)^2+4(k-\mu)$, we get 
\begin{equation}\label{delta_cond}
\Delta=\left(\lambda-\frac{k(k-\lambda-1)}{v-k-1}\right)^2+4\left(k-\frac{k(k-\lambda-1)}{v-k-1}\right).
\end{equation}

\begin{lemma}\label{k_16p2}
If a non-trivial $(v,k,\lambda, \mu)$-PDS exists in an Abelian group with $v=8p^3$, $p\geq5$ a prime, $\Delta=16p^2$ and $k\le \frac{v}{2}$, then we have $k\le 4p^2+3p-\frac{1}{2}$.
\end{lemma}

\proof Setting $\Delta=16p^2$ and $v=8p^3$ in Equation (\ref{delta_cond}),  and solving the obtained quadratic equation for $\lambda$, we get
$$\lambda=\frac{-3k(8p^3-1)+k^2(8p^3+1)\pm 2\sqrt{(1+k-8p^3)^2 \,(4p^2(8p^3-1)^2+k(1+k-8p^3)8p^3})}{(8p^3-1)^2}.$$

As $\lambda$ is an integer, the discriminant must be nonnegative, hence
\begin{equation}\label{ineq_1}
4p^2(8p^3-1)^2+k(1+k-8p^3)8p^3\ge 0.
\end{equation}

Solving the quadratic equation $4p^2(8p^3-1)^2+k(1+k-8p^3)8p^3= 0$ for $k$ yields
$$k=\frac{-p+8p^4\pm \sqrt{64p^8-128p^7-16p^5+32p^4+p^2-2p}}{2p}.$$

As the largest root is greater than $4p^3$, and $k\leq v/2$, in order for Inequality (\ref{ineq_1}) to hold,  we have 
$$k \le \frac{-p+8p^4- \sqrt{64p^8-128p^7-16p^5+32p^4+p^2-2p}}{2p}.$$
We estimate $64p^8-128p^7-16p^5+32p^4+p^2-2p=(8p^4-8p^3-6p^2)^2+32p^6-112p^5-4p^4+p^2-2p\ge (8p^4-8p^3-6p^2)^2$ for $p \ge 5$. 
Hence $$k \le \frac{-p+8p^4-\sqrt{(8p^4-8p^3-6p^2)^2}}{2p}=4p^2+3p-\frac{1}{2}.$$
\qed

\begin{lemma}\label{k_4p2}
If a non-trivial $(v,k,\lambda, \mu)$-PDS exists in an Abelian group with $v=8p^3$, $p\geq5$ a prime, $\Delta=4p^2$ and $k\le \frac{v}{2}$, then we have $k\le p^2+\frac{p}{4}-\frac{1}{2}$.
\end{lemma}

\proof This proof is very similar to that of the preceding lemma. Setting $\Delta=4p^2$ and $v=8p^3$ in Equation (\ref{delta_cond}),  and solving the obtained quadratic equation for $\lambda$, we now obtain 
\begin{equation}\label{ineq_2}
p^2(8p^3-1)^2+k(1+k-8p^3)8p^3\ge 0.
\end{equation}

Solving the quadratic equation $p^2(8p^3-1)^2+k(1+k-8p^3)8p^3= 0$ for $k$ yields
$$k=\frac{-2p+16p^4\pm \sqrt{256p^8-128p^7-64p^5+32p^4+4p^2-2p}}{4p}.$$

Again as the largest root is greater than $4p^3$, and $k\leq v/2$ is assumed, in order for Inequality (\ref{ineq_2}) to be true,  we have 
$$k \le \frac{-2p+16p^4- \sqrt{256p^8-128p^7-64p^5+32p^4+4p^2-2p}}{4p}.$$
We estimate $256p^8-128p^7-64p^5+32p^4+4p^2-2p=(16p^4-4p^3-p^2)^2+16p^6-72p^5+31p^4+4p^2-2p\ge (16p^4-4p^3-p^2)^2$  for $p \ge 5$. 
Hence $$k \le \frac{-2p+16p^4-\sqrt{(16p^4-4p^3-p^2)^2}}{4p}=p^2+\frac{p}{4}-\frac{1}{2}.$$
\qed

\subsection{Restrictions on $\mu$}

From Lemma \ref{Delta} we know that we may assume that $\Delta$ equals either $4p^2$ or $16p^2$. In this subsection we will show that $\Delta=4p^2$ cannot occur, and that either $\mu=2p+2$ or  $\mu=2p-2$ when $\Delta=16p^2$.

\begin{proposition}
The case $\Delta=4p^2$ cannot occur.
\end{proposition}

\proof If $\Delta=4p^2$, by Proposition \ref{Ma97}-(e), we have $(2k-\beta)^2\equiv 0 \pmod {4p^2}$. Thus we can write $k=px+\frac{\beta}{2}$ with $x$ an integer.

By substituting  $v=8p^3$,  $\lambda=\mu+\beta$, and $k=px+\frac{\beta}{2}$, Equation (\ref{srg_cond_1}) becomes 
\begin{equation}\label{srg_8pcube_1}
p^2x^2-\frac{\beta^2}{4}-px-\frac{\beta}{2}-8p^3\mu+\mu=0.
\end{equation}

As $\Delta=4p^2=(\lambda-\mu)^2+4(k-\mu)$=$\beta^2+4(px+\frac{\beta}{2}-\mu)$, it follows that
\begin{equation}\label{beta_cond_1}
\beta^2+2\beta=4p^2-4px+4\mu.
\end{equation}

Combining  Equations (\ref{srg_8pcube_1}) and (\ref{beta_cond_1}), we have obtained
\begin{equation}\label{mu_cond_1}
\mu=\frac{x^2-1}{8p}.
\end{equation}

Since $\mu$ is an integer,  $x+1$ and $x-1$ have the same parity, and $p$ is a prime number $\ge 5$, it follows that $2p\,|\, x+1$ or $2p \,|\, x-1$. So we have $x=2tp-1$ or $x=2tp+1$, where $t$ is an integer. Thus $k=px+\frac{\beta}{2}=2tp^2\pm p+\frac{\beta}{2}$. Since $k \le p^2+\frac{p}{4}-\frac{1}{2}$ by Lemma \ref{k_4p2}  and  $-\sqrt{\Delta }<\beta <\sqrt{\Delta}-2$ by Proposition \ref{Ma97}, we have $t=0$. But 
when $t=0$, we have  $x=\pm 1$ and $\mu=0$, implying a trivial PDS. \qed

\medskip

\begin{lemma}
If $\Delta=16p^2$, then either $\mu=2p+2$ or  $\mu=2p-2$.
\end{lemma}

 \proof If $\Delta=16p^2$, by Proposition \ref{Ma97}-(e), we have $(2k-\beta)^2\equiv 0 \pmod {16p^2}$. Thus we can write $k=2px+\frac{\beta}{2}$ with $x$ an integer.

By substituting  $v=8p^3$,  $\lambda=\mu+\beta$, and $k=2px+\frac{\beta}{2}$, Equation (\ref{srg_cond_1}) becomes 
\begin{equation}\label{srg_8pcube_2}
4p^2x^2-\frac{\beta^2}{4}-2px-\frac{\beta}{2}-8p^3\mu+\mu=0.
\end{equation}

As $\Delta=16p^2=(\lambda-\mu)^2+4(k-\mu)$=$\beta^2+4(2px+\frac{\beta}{2}-\mu)$, it follows that
\begin{equation}\label{beta_cond_2}
\beta^2+2\beta=16p^2-8px+4\mu.
\end{equation}

Combining  Equations (\ref{srg_8pcube_2}) and (\ref{beta_cond_2}), we have obtained
\begin{equation}\label{mu_cond_2}
\mu=\frac{x^2-1}{2p}.
\end{equation}

Since $\mu$ is an integer,  $x+1$ and $x-1$ have the same parity and $p$ is a prime number $\ge 5$, it follows that $2p\,|\, x+1$ or $2p \,|\, x-1$. So we can assume that $x=2tp-1$ or $x=2tp+1$, where $t$ is an integer. If $t= 0$  then $x=\pm 1$ and $\mu=0$, and we obtain a trivial PDS.

Thus $k=2px+\frac{\beta}{2}=4tp^2\pm 2p+\frac{\beta}{2}$. Since $k \le 4p^2+3p-\frac{1}{2}$ by Lemma \ref{k_16p2}  and  $-\sqrt{\Delta }<\beta <\sqrt{\Delta}-2$ by Proposition \ref{Ma97}, we have $t=1$. Hence either $x=2p+1$ and $\mu=2p+2$; or $x=2p-1$ and  $\mu=2p-2$. \qed

\section{Non-existence of PDS in Abelian groups of order $8p^3$ with $\mu=2p+2$}
When $x=2p+1$, we have  $\mu=\frac{x^2-1}{2p}=2p+2$ and $k=2px+\frac{\beta}{2}=4p^2+2p+\frac{\beta}{2}$. Also Equation (\ref{beta_cond_2}) becomes
\begin{equation}\label{beta_new_cond}
\beta^2+2\beta-8=0
\end{equation}
Thus $\beta=-4$ or $\beta=2$. If $\beta=-4$, then $k=4p^2+2p-2$, $\mu=2p+2$,  and $\lambda=2p-2$;  if $\beta=2$, then $k=4p^2+2p+1$, $\mu=2p+2$,  and $\lambda=2p+4$.

\begin{theorem}\label{8pcube_case1}
There does not exist a regular ($8p^3$, $4p^2+2p-2$, $2p-2$, $2p+2$)--PDS in an Abelian group, where $p\ge 5 $ is a prime.
\end{theorem}

\proof Let $G$ be an Abelian group of order $8p^3$. Assume by way of contradiction that $D$ is a ($8p^3$, $4p^2+2p-2$, $2p-2$, $2p+2$)--PDS in $G$. By Proposition \ref{Ma2}, we know that $G\cong\mathbb{Z}_2^3 \times\mathbb{Z}_p^3$. 
Let $g_1$, $g_2$, $\cdots$, $g_{p^3-1}$ be all elements of order $p$ in $G$, and let $\cB_{g_i}=\{ag_i\,|\,o(a)=1 \;\mbox{or}\; 2, \; ag_i\in D\}$, and $B_i=|\cB_{g_i}|$, $i=1,2, \cdots, p^3-1$. That is, $B_i$ equals the number of elements in $D$ whose $(p+1)$-th power equals $g_i$.  We next prove the simple yet important observation that the LMT implies that $|\cB_{g_i}|=|\cB_{g_i^s}|$  whenever $1 \le s \le p-1$.

\medskip

\noindent {\bf Observation (O)}   {\it It holds that $|\cB_{g_i}|=|\cB_{g_i^s}|$ for $1 \le s \le p-1$.}
\smallskip

\proof If $ag_i \in \cB_{g_i}$, that is, if $ag_i \in D$ with $o(a)=1$ or 2, then, by the LMT, $(ag_i)^s=ag_i^s\in D$ if s is odd;
or $(ag_i)^{s+p}=ag_i^s\in D$ if s is even. Thus $|\cB_{g_i}|\le |\cB_{g_i^s}|$.

On the other hand, since $\gcd(s,p)=1$, we can find integers $r$ and $t$ such that $rs+tp=1$. It is clear that $\gcd(r,p)$=$\gcd(r+p,p)$=1. If $ag_i^s\in \cB_{g_i^s}$, that is, if $ag_i^s \in D$ with $o(a)=1$ or 2, then, by the LMT, $(ag_i^s)^r=ag_i \in D$ if $r$ is odd, and $(ag_i^s)^{r+p}=ag_i\in D$ if $r$ is even. Thus $|\cB_{g_i^s}|\le | \cB_{g_i}|$.  

Hence $|\cB_{g_i}|=|\cB_{g_i^s}|$ for $1\le s \le p-1$. \qed

\medskip

Let $N$ be the Sylow 2-subgroup of $G$. Using Proposition \ref{Ma1}, we know that $D$ contains either 0 or 4 elements of order 2. First assume that $D$ contains no elements of order 2. We see that $\Sigma_i B_i=4p^2+2p-2$ and $\Sigma_iB_i(B_i-1)=14p+14$, as all 7 elements of order 2 are not in $D$, thus each  of them has exactly $\mu=2p+2$ difference representations.\

According to observation (O), we can assume, by relabeling the $g_i$ if necessary, that $C_j:=B_{(j-1)(p-1)+1}=B_{(j-1)(p-1)+2}=\cdots=B_{(j-1)(p-1)+(p-1)}$ for $j=1,2,\cdots, p^2+p+1$, and $C_1\ge C_2\ge \cdots \ge C_{p^2+p+1}$. We now obtain 
\begin{eqnarray}\label{2p-2_Case1}
\Sigma_jC_j&=&\frac{4p^2+2p-2}{p-1}=(4p+6)+\frac{4}{p-1},\\
 \Sigma_j C_j^2&=&\frac{4p^2+16p+12}{p-1}.\label{2p-2_Case1-II}
\end{eqnarray}

Since $\frac{4}{p-1}$ is an integer only when $p-1$=1, 2, or 4, that is, when $p$=2, 3, or 5, Equation (\ref{2p-2_Case1}) has no integer solutions when $p\ge 7$ and $p$ is a prime number. 

When $p=5$,  Equations (\ref{2p-2_Case1}) and (\ref {2p-2_Case1-II}) become
$\Sigma_jC_j=27$ and 
 $\Sigma_j C_j^2=48$. Thus $\Sigma_{j} C_j (C_j-1)=21$, which contradicts with the fact that $C_j(C_j-1)$ is always even.

\medskip

 Secondly assume that $D$ contains $4$ elements of order 2. It follows that $\Sigma_i B_i+4=4p^2+2p-2$. By counting the number of difference representations of elements of order 2 in $DD^{(-1)}$, we obtain that  $\Sigma_iB_i(B_i-1)+4\cdot 3=4(2p-2)+3(2p+2)$. Using similar labelling as above, we now obtain 
\begin{equation}\label{2p-2_Case2}
\Sigma_{j=1}^{p^2+p+1}C_j=4p+6 \quad \mbox {and} \quad \Sigma_{j=1}^{p^2+p+1} C_j^2=4p+20.
\end{equation}

Since $C_j$ is a non-negative integer, and $0^2-0=0$, $1^2-1=0$, $2^2-2=2$, $3^2-3=6$, $4^2-4=12$, $5^2-5=20$, $\Sigma_{j=1}^{p^2+p+1} (C_j^2-C_j)=14$, the system of Equations (\ref{2p-2_Case2}) only has the following nonnegative integer solutions, listed as (decreasing) $p^2+p+1$ tuples:

$(4, 2,\, \underbrace{1, 1, \cdots, 1}_{4p}\,, 0, 0, \cdots, 0)$,\\

 $(3, 3, 2,\, \underbrace{1, 1, \cdots, 1}_{4p-2}\,, \,0, 0, \cdots, 0)$,\\

$(3, 2, 2, 2, 2,\,\underbrace{ 1, 1, \cdots, 1}_{4p-5}\,, 0, 0, \cdots, 0)$,\\

$(2, 2, 2, 2, 2, 2, 2,\,\underbrace{ 1, 1 , \cdots, 1}_{4p-8}\,, 0, 0, \cdots, 0)$.\\

Recall that $N$ is the unique subgroup isomorphic to $\mathbb{Z}_2^3$ in $G$. Let $P_1,\hdots,P_{p^2+p+1}$ be the $p^2+p+1$ subgroups of $G$ isomorphic to $\mathbb{Z}_p$, and let $L_1,\hdots,L_{p^2+p+1}$ be the $p^2+p+1$ subgroups of $G$ isomorphic to $\mathbb{Z}_p^2$. Now consider the incidence structure $\mathcal{P}$ with points the subgroups $P_i\times N$, $i=1,\hdots,p^2+p+1$,  of $G$,  with blocks the subgroups $L_i\times N$, $i=1,\hdots, p^2+p+1$, of $G$, and with containment as incidence. Then it is easily seen that $\mathcal{P}$ is a $2-(p^2+p+1,p+1,1)$ design, or equivalently, the unique projective plane of order $p$. We next assign a weight to each point of $\mathcal{P}$ in the following way: if point $p$ corresponds to  subgroup $P_i\times N$ then the weight of $p$ is $\frac{1}{p-1}|((P_i\times N)\setminus N)\cap D|$. In this way the weights of the $p^2+p+1$ points of $\mathcal{P}$ correspond to the $p^2+p+1$ values $C_{1},C_{2},\hdots, C_{p^2+p+1}$, that is, $1/(p-1)$ times the number of elements of order $p$ or $2p$ from $ D$ in the subgroup underlying the given point. Without loss of generality we may assume the labeling is such that  point $P_i\times N$ has weight $C_i$. The weight of a block will simply be the sum of the weights of the points in that block.
\medskip

We next count how many elements of order $p$ or $2p$ from $D$ a specific subgroup of the form $L_i\times N$ can contain. Assume that $|(L_i\times N) \cap D|=m$. Let $ag$ and $bh$ be two distinct elements from $ D$, with $a^2=b^2=g^p=h^p=e$. Then $agh^{-1}b^{-1}$ belongs to $L_i\times N$ if and only if $gh^{-1}\in L_i$. It is easy to see that if $g\in L_i$ there are $m-1$ possibilities for $bh$ such that $gh^{-1}\in L_i$. When $g\not\in L_i$ there are several cases to discuss (recall that, by the LMT, when an element $bh$ belongs to $D$, so do all elements $bh^l$ for $1\leq l\leq p-1$) :
\begin{itemize}
\item if $bh\in L_i\times N$ then obviously $agh^{-1}b^{-1}$ does not belong to $L_i\times N$;
\item if $bh=ag^l$ for some $l\in\{1,\hdots,p-1\}$ then clearly $agh^{-1}b^{-1}$ cannot be a nonidentity element of $L_i\times N$;
\item if $bh=bg^l$ for some $l\in\{1,\hdots,p-1\}$ and $b\neq a$ then $agh^{-1}b^{-1}$ will be in $L_i\times N$ if and only if $l=1$;
\item if $h\neq g^l$ for any $l\in\{1,\hdots,p-1\}$ and $h\notin L_i$, then it is easy to see there is a unique $r\in\{1,\hdots,p-1\}$ such that $gh^{-r}\in L_i$, and hence such that $agh^{-r}b^{-1}$ belongs to $L_i\times N$.
\end{itemize}
Combining the above observations yields that when  $g\notin L_i$ there are $\frac{| D|-m-(p-1)}{p-1}$ possibilities for $bh$ such that $agh^{-1}b^{-1}\in L_i\times N$.
\medskip

 Counting the number of differences of elements of $ D$ that are in $L_i\times N$ in two ways, we obtain 

\begin{equation}\label{key}
m(m-1)+(k-m)(\frac{k-m-(p-1)}{p-1})=\lambda m +\mu(8p^2-1-m),
\end{equation}
where $(k,\lambda,\mu)=(4p^2+2p-2, 2p-2, 2p+2)$. This yields that $m=2(p+1)$ or $2(3p-1)$.
\medskip

Now define $m':=\frac{1}{p-1}|((L_i\times N)\setminus N) \cap D|$. We obtain $m'=\frac{2(p+1)-4}{p-1}=2$ or $m'=\frac{2(3p-1)-4}{p-1}=6$ since $D$ contains 4 elements of order 2.
\medskip

We now note that the values $m'$ must be the weights of the blocks of $\mathcal{P}$, and that in both cases these weights are even. We first show that no value $C_{i}$ can be odd. Assume by way of contradiction that $C_{i}$ is odd for some $i$. Let the weight of the $p+1$ blocks that contain $P_i\times N$ be $n_{1},\hdots, n_{p+1}$ respectively. Then 

$$\sum_{j=1}^{p^2+p+1}C_{j}=C_{i}+\sum_{t=1}^{p+1}(n_t-C_{i}).$$

As $n_t$ is even for all $t$ (the $n_t$ are $m'$ values), this implies that $\sum_{j=1}^{p^2+p+1}C_{j}$ is odd. This contradicts with the fact that $\sum_{j=1}^{p^2+p+1}C_{j}=4p+6$.
\medskip

Since all the solutions to the system of Equations (\ref{2p-2_Case2}) contain at least one odd $C_j$, it follows that no nontrivial ($8p^3$, $4p^2+2p-2$, $2p-2$, $2p+2$)--PDS exists in an Abelian group.  \qed

\bigskip

\begin{theorem}\label{case2}
There does not exist a ($8p^3$, $4p^2+2p+1$, $2p+4$, $2p+2$)--PDS  in an Abelian group, where $p\ge 5 $ is a prime. 
\end{theorem}

\proof This case is dealt with in a very similar way. We will only provide a sketch of the proof. Assume by way of contradiction $D$ is a ($8p^3$, $4p^2+2p+1$, $2p+4$, $2p+2$)-PDS in an Abelian group $G$.

As before $G\cong\mathbb{Z}_2^3 \times\mathbb{Z}_p^3$, and using Proposition \ref{Ma1} we obtain that $D$ contains either $3$ or $7$ elements of order $2$.
 If $D$ contains $3$ elements of order $2$ we obtain 
\begin{eqnarray}\label{2p+4_Case1}
\Sigma_jC_j&=&\frac{4p^2+2p-2}{p-1}=(4p+6)+\frac{4}{p-1},\\
 \Sigma_j C_j^2&=&\frac{4p^2+16p+12}{p-1}.\label{8pcube_Case1-II}
\end{eqnarray}
 which is the same as the system of Equations (\ref{2p-2_Case1}), (\ref{2p-2_Case1-II}). Hence this case cannot occur.\

 If $D$ contains $7$ elements of order $2$ we obtain 
\begin{equation}\label{2p+4_Case2}
\Sigma_jC_j=4p+6 \quad \mbox {and} \quad \Sigma_j C_j^2=4p+20,
\end{equation}
which is the same as the system of equations in (\ref{2p-2_Case2}),  and thus  has the same set of solutions.

With similar notation as in the previous theorem, and using the same counting argument for $(k,\lambda,\mu)$=($4p^2+2p+1$, $2p+4$, $2p+2$), one obtains $m=2p+5$ or $m=6p+1$.
\medskip

Now define $m':=\frac{1}{p-1}|((L_i\times N)\setminus N) \cap D|$. We obtain $m'=\frac{2p+5-7}{p-1}=2$ or $m'=\frac{6p+1-7}{p-1}=6$ since $D$ contains 7 elements of order 2.
\medskip

As before the weights of all blocks of $\mathcal{P}$ must be even, and the proof can be finished in the same way as in the ($8p^3$, $4p^2+2p-2$, $2p-2$, $2p+2$) case.\qed

\bigskip

\section{Non-existence of PDS in Abelian groups of order $8p^3$ with $\mu=2p-2$}
When $x=2p-1$, we have  $\mu=\frac{x^2-1}{2p}=2p-2$ and $k=2px+\frac{\beta}{2}=4p^2-2p+\frac{\beta}{2}$. Also Equation (\ref{beta_cond_2}) becomes
\begin{equation}\label{beta_new_cond}
\beta^2+2\beta-8(2p-1)=0
\end{equation}
Thus $\beta=-1\pm \sqrt{16p-7}$. It is easy to observe that when $\sqrt{16p-7}$ is an integer it is congruent to 3 or 5 modulo 8:
\begin{itemize}
\item
If $\sqrt{16p-7}=8y+1$, we have $p=4y^2+y+\frac{1}{2}$, contradicting with the assumption that $p$ is an integer;
\item
If $\sqrt{16p-7}=8y+3$, we have $p=4y^2+3y+1$;
\item
If $\sqrt{16p-7}=8y+5$, we have $p=4y^2+5y+2$; 
\item
If $\sqrt{16p-7}=8y+7$, we have $p=4y^2+7y+\frac{7}{2}$, contradicting with the assumption that $p$ is an integer.

\end{itemize}

\subsection  {The case $p=4y^2+3y+1$ }

In this subsection, we let  $p=4y^2+3y+1$ and $p\ge 5$ be a prime number.  It is easy to see that $y\ge 2$,  and $\beta=-1\pm (8y+3)$=$-8y-4$ or $8y+2$ by Equation (\ref{beta_new_cond}).
\begin{theorem}\label{case3}
Let $p=4y^2+3y+1$ and $p\ge 5$ be a prime number. Then no non-trivial regular  partial difference sets exist with $\mu=2p-2$, $\beta=-8y-4$ or $8y+2$, and $k=4p^2-2p+\frac{\beta}{2}$ in Abelian groups of order $8p^3$.
\end{theorem}

 \proof  Let $G$ be an Abelian group of order $8p^3$. We will prove this theorem in two parts based on the $\beta$ values:

\begin{itemize}
\item[(i)] Let $\beta=-8y-4$. Assume on the contrary that $D$ is a non-trivial regular $(8p^3,\,4p^2-2p-4y-2, \, 2p-8y-6, \, 2p-2)$-PDS in $G$. Assume that $|D \cap Z_2^3|=a$, where $0\le a \le 7$.  Using the notation from Section 3, we have 
$$\sum B_i=k-a=4p^2-2p-4y-2-a.$$
It follows that $$\sum C_i=\frac{4p^2-2p-4y-2-a}{p-1}=4p+2-\frac{4y+a}{p-1}.$$
Since $y\ge 2$ and $0\le a \le 7$, and $p=4y^2+3y+1$, it is easy to check that
$$0< 4y+a \le 4y+7 < p-1=4y^2+3y,$$
thus $\frac{4y+a}{p-1}$ is not an integer, and hence no such $D$ exists.

\item[(ii)] Let $\beta=8y+2$. Assume on the contrary that $D$ is a non-trivial  regular $(8p^3, \, 4p^2-2p+4y+1, \, 2p+8y, \, 2p-2)$  PDS in $G$. Assume that $|D \cap Z_2^3|=a$, where $0\le a \le 7$.  Using the  notation from Section 3, we have 
$$\sum B_i=4p^2-2p+4y+1-a.$$
It follows that $$\sum C_i=\frac{4p^2-2p+4y+1-a}{p-1}=4p+2+\frac{4y+3-a}{p-1}.$$
Since $y\ge 2$ and $0\le a \le 7$, and $p=4y^2+3y+1$, it is easy to check that
$$0< 4y+3-a \le 4y+3 < p-1=4y^2+3y,$$
thus $\frac{4y+3-a}{p-1}$ is not an integer, and hence no such $D$ exists.
\end{itemize}
\qed

\subsection  {The case $p=4y^2+5y+2$}
In this subsection, we let  $p=4y^2+5y+2\ge 5$ be a prime number. It is easy to see that  $\beta=-1\pm (8y+5)$=$-8y-6$ or $8y+4$ by Equation (\ref{beta_new_cond}).
We note that the first part of the following theorem is slightly more subtle than the proof of Theorem \ref{case3} as one needs to invoke Proposition  \ref{Ma1} in order to exclude the case of $5$ elements of order $2$ when $p=11$.

\begin{theorem}\label{case4}
Let $p=4y^2+5y+2$ and $p\ge 5$ be a prime number. Then no non-trivial  regular partial difference sets exist with $\mu=2p-2$, $\beta=-8y-6$ or $8y+4$, and $k=4p^2-2p+\frac{\beta}{2}$ in Abelian groups of order $8p^3$.
\end{theorem}

 \proof  Let $G$ be an Abelian group of order $8p^3$. We will prove this theorem in two parts based on the $\beta$ values:

\begin{itemize}
\item[(i)] Let $\beta=-8y-6$.  Assume on the contrary that $D$ is a non-trivial  regular $(8p^3, \, 4p^2-2p-4y-3, \, 2p-8y-8, \, 2p-2)$ regular PDS in $G$.

Let $N=Z_2^3$ be the Sylow 2-group of $G$. By Proposition \ref{Ma1}, we know that $D$ contains either 3 or 7 elements of order 2. Thus $$\sum{B_i}=4p^2-2p-4y-6 \quad \mbox{ or} \quad  \sum{B_i}=4p^2-2p-4y-10.$$ Hence 
$$\sum{C_i}=4p+2-\frac{4y+4}{p-1} \quad \mbox {or} \quad  \sum{C_i}=4p+2-\frac{4y+8}{p-1}.$$

Clearly $0<4y+4 <p-1=4y^2+5y+1$ when $y \ge 1$, so $\frac{4y+4}{p-1}$ is not an integer. Also  $0<4y+8<p-1=4y^2+5y+1$ when $y\ge 2$,  and $\frac{4y+8}{4y^2+5y+1}=\frac{12}{10}$ when $y=1$, so $\frac{4y+8}{4y^2+5y+1}$ is not an integer for any $y\ge 1$.  Hence no such $D$ exists.

\item[(ii)]  Let $\beta=8y+4$. Assume on the contrary that $D$ is a non-trivial  regular $(8p^3, \, 4p^2-2p+4y+2, \, 2p+8y+2, \, 2p-2)$ regular PDS in $G$.  Assume that $|D \cap Z_2^3|=a$, where $0\le a \le 7$. Then 
$$\sum B_i=4p^2-2p+4y+2-a.$$
It follows that $$\sum C_i=\frac{4p^2-2p+4y+2-a}{p-1}=4p+2+\frac{4y+4-a}{p-1}.$$
Since $y\ge 1$ and $0\le a \le 7$, and $p=4y^2+5y+2$, it is easy to check that
$$0< 4y+4-a \le 4y+4 < p-1=4y^2+5y,$$
thus $\frac{4y+4-a}{p-1}$ is not an integer. Hence no such $D$ exists.
\end{itemize}

\qed

\section{Conclusion}

Combining all results from this paper and \cite{SDWZW_216} we have obtained the following theorem.

\begin{theorem}
All regular PDS in Abelian groups of order $8p^3$, $p$ an odd prime, are trivial.
\end{theorem}
\proof When $p=3$ this is the main result of \cite{SDWZW_216}. When $p\geq 5$ this follows by combining Theorems \ref{8pcube_case1}, \ref{case2} ,\ref{case3}, \ref{case4}. \qed

As mentioned in the introduction, this nonexistence result seems to provide further evidence that regular PDS in Abelian groups are rare. Nevertheless, classifying or completely characterizing regular PDS in Abelian groups seems to be completely out of reach at this point.

\end{document}